\documentclass[12pt,a4paper,reqno]{amsart}

\usepackage{xcolor}
\usepackage[all,2cell,arrow,color]{xy}

\usepackage{amssymb}
\usepackage{mathtools}
\usepackage{tensor}
\usepackage{graphicx} 
\usepackage{rotating}
\usepackage{pdflscape}
\definecolor{darkgreen}{rgb}{0,0.45,0} 
\usepackage[pagebackref,colorlinks,citecolor=darkgreen,linkcolor=darkgreen]{hyperref}
\usepackage{enumerate,xspace}

\def\beq#1{\begin{equation}\label{#1}}
\def\eeq{\end{equation}}

\newcommand{\bb}{\ensuremath{\mathbb B}\xspace}
\newcommand{\bc}{\ensuremath{\mathbb C}\xspace}
\newcommand{\cv}{\ensuremath{\mathcal V}\xspace}
\newcommand{\cb}{\ensuremath{\mathcal B}\xspace}

\newcommand{\Set}{\ensuremath{\mathsf{Cat}}\xspace}
\newcommand{\Cat}{\ensuremath{\mathsf{Cat}}\xspace}
\newcommand{\Mat}{\ensuremath{\mathsf{Mat}}\xspace}
\newcommand{\Span}{\ensuremath{\mathsf{Span}}\xspace}

 \newtheorem{proposition}{Proposition}[section]

  \newtheorem{theorem}[proposition]{Theorem}

  \theoremstyle{definition}
  \newtheorem{definition}[proposition]{Definition}
  \newtheorem{example}[proposition]{Example}

  \theoremstyle{remark}

  \newcounter{c}
  \renewcommand{\[}{\setcounter{c}{1}$$}
  \newcommand{\etyk}[1]{\vspace{-7.4mm}$$\begin{equation}\Label{#1}
  \addtocounter{c}{1}}
  \renewcommand{\]}{\ifnum \value{c}=1 $$\else \end{equation}\fi}
\title{A note on warpings of monoidal structures}
\author{Dimitri Chikhladze} 
\address{}
\email{d.chikhladze@gmail.com}
\date{}

\begin{document}

\maketitle
\begin{abstract}
In \cite{LS14} the analogy between the Kleisli construction and the construction of ``warping a skew monoidale category'' in the sense of \cite{LS12} was outlined. In this note we present the same work in a slightly more formal way.
\end{abstract}

\section{Warpings}\label{Warp}
Fix a bicategory $\bb$ as a large ambient. When we talk of objects, morphisms, 2-cells, monads and other concepts we mean those within $\bb$. As usual we omit the bicategory structural isomorphisms.

\begin{definition}
A \textbf{warping} from a monad $(X, B)$ to an endomorphism $A : X \rightarrow X$ consists of 2-cells  $t : ABA \Rightarrow AB$ and $k : 1_X \rightarrow AB$ satisfying the following equations

\[
\xymatrix@R=1em{
&ABBA\ar@{=>}[r]^{ApA}&ABA\ar@{=>}[dr]^{t}&\\
ABABA \ar@{=>}[rd]_{ABt}\ar@{=>}[ur]^{tBA}&&&AB\\
&ABAB\ar@{=>}[r]_{tB}&ABB\ar@{=>}[ur]_{Ap}&\\
}
\]

\[\xymatrix{ABAB \ar@{=>}[r]^{tB}&ABB\ar@{=>}[d]^{Ap}\\
AB\ar@{=>}[u]^{ABk}\ar@{=>}[r]_{1_{AB}}&AB}
\qquad\qquad
\xymatrix{&ABA\ar@{=>}[rd]^{t}&\\
A\ar@{=>}[ur]^{kA}\ar@{=>}[rr]_{Ae}&&AB
}
\]

\end{definition}

A warping is precisely the structure which is needed to define a monad structure on the composite $AB$ which is suitably compatible with the monad $(X, B)$. 

\begin{theorem}\label{warptheo}
Given a monad $(X, B)$ and an endomorphism $A : X \rightarrow X$, there is a one-to-one correspondence between warping structures $(t, k)$ and those monad structures on $AB$ for which $(Ap)(eB): B \rightarrow AB$ is a monoid map; a warping is determined from a monoid structure on $AB$ by
\begin{itemize}
\item $t = p^{AB}(ABAe)$
\item $k = e^{AB}$
\end{itemize}
\noindent and conversely, a monad structure on $AB$ is determined by a warping by
\begin{itemize}
\item $p^{AB} = (Ap)(tB)$
\item $e^{AB} = k$.
\end{itemize}
\end{theorem}

Recall from \cite{LS02} that a wreath consists of a monad $(X, B)$, an endofunctor $A : X \rightarrow X$, and 2-cells $d : BA \Rightarrow AB$, $q : AAB \Rightarrow AB$, $j : 1_X \Rightarrow AB$ satisfying the following axioms of which the first pair expresses compatibility of $d$ with the monad structure of $B$, the second pair expresses compatibility of $d$ with $q$ and $j$, while the last set in a sense represents associativity and unitivity conditions for $q$ and $j$:

\[\xymatrix{
BBA\ar@{=>}[r]^{Bd}\ar@{=>}[d]_{pA}&BAB\ar@{=>}[r]^{dB}&ABB\ar@{=>}[d]^{Ap}\\
BA \ar@{=>}[rr]_{d}&&AB\\
}
\qquad
\xymatrix{&A\ar@{=>}[ld]_{eA}\ar@{=>}[rd]^{Ae}&\\
BA\ar@{=>}[rr]_{d}&&BA}
\]

\[
\xymatrix{BAA \ar@{=>}[rd]^{Bq} \ar@{=>}[r]^{dA}&ABA \ar@{=>}[r]^{Ad}&AAB\ar@{=>}[r]^{qB}&ABB \ar@{=>}[d]^{Ap}\\
&BAB\ar@{=>}[r]^{dB}&ABB\ar@{=>}[r]^{Ap}&AB}
\qquad
\xymatrix{B \ar@{=>}[rr]^{jB}\ar@{=>}[d]_{Bj}&&ABB\ar@{=>}[d]^{Ap}\\
BAB\ar@{=>}[r]_{dB}&ABB\ar@{=>}[r]_{Ap}&AB}
\]

\[
\xymatrix{AAA \ar@{=>}[r]^{Aq}\ar@{=>}[d]_{qA}&AAB\ar@{=>}[r]^{qB}&ABB\ar@{=>}[dd]^{Ap}\\
ABA\ar@{=>}[d]_{Ad}&&\\
AAB\ar@{=>}[r]_{qB}&ABB\ar@{=>}[r]_{Ap}&AB\\}
\qquad
\xymatrix{A\ar@{=>}[r]^{Aq}\ar@{=>}[rdd]_{Ae}&AAB\ar@{=>}[d]^{qB}&ABA\ar@{=>}[l]_{Ad}&A\ar@{=>}[l]_{j}\ar@{=>}[lldd]^{Ae}\\
&ABB\ar@{=>}[d]^{Ap}&&\\
&AB&&}
\]

\noindent  In fact, a wreath is a monad in a certain bicategory $\mathrm{EM}(\bb)$ (see \cite{LS02}). Recall also, that a wreath between $(X, B)$ and $A$ gives rise to a monad $(X, AB)$ with the monad unit $j : 1 \Rightarrow AB$ and the monad multiplication $(pA)(qp)(AdB) : ABAB \Rightarrow AB$.

\begin{theorem}\label{warptowreath}
There is a one-to-one correspondence between warpings $(t, k)$ from a monad $(X, B)$ to an endofunctor $A$ and wreaths $(d, q, j)$ between the same two; a wreath is determined from a warping by
\begin{itemize}
\item $d = t(Ap A)(kBA)$
\item $q = t(AeA)$
\item $j = k$
\end{itemize}
\noindent and conversely, a warping is determined from a wreath by 
\begin{itemize}
\item $t = (Ap)(qB)(Ad)$
\item $k = j$.
\end{itemize}
\end{theorem}

\begin{example}\label{exwarp}
Take $\bb$ to be the bicategory of spans $\Span$. Its objects are sets, a morphism $M: X \rightarrow Y$ in it can be identified with a collection $M(x, y)$ of sets indexed by pairs $(x, y) \in X\times Y$, while its 2-cells $M \Rightarrow M'$ can be identified with families of maps $M(x, y) \rightarrow M'(x, y)$. The composition is expressed by the usual matrix multiplication formula
\[MN(x, y) = \coprod_zM(z, y)\times N(x, z).\]

\noindent The identity morphism $1_X$ is the collection of sets indexed by the elements of $X\times X$ consisting of the one-element set on the diagonal and the empty set everywhere else. There is a pseudofunctor $\Set^{\mathrm{op}} \rightarrow \Span$ which is identical on objects, and sends a map $F : X \rightarrow Y$ to the span $F^\ast$ for which $F^\ast(y, x)$ is the one-element set whenever $Fx = y$ and the empty set otherwise. $F^\ast$ has a dually defined left adjoint $F_\ast$. We have
\[F^\ast M(z, x) = M(z, Fx)\]
\[NF^\ast(x, z) = N(Fx, z)\]

\noindent from which it follows that, a 2-cell  $MF^\ast \rightarrow F^\ast M$ amounts to a collection of maps $M(x, y) \rightarrow M(Fx, Fy)$.

It is a well-known simple fact that a monad $(X, B)$ in $\Span$ is a category whose set of object is $X$, and whose homsets are $B(x, y)$. We use $B$ to stand for this category as well.

A wreath $(d, q, j)$ from $(X, B)$ to $T^\ast$ within $\Span$ is a monad on the category $B$. The 2-cell $d$ in this case amounts to a family $B(x, y) \rightarrow B(Tx, Ty)$, which gives the extension of the function $T$ to a functor on the category $B$. The 2-cells $q$ and $j$ amount to specifying natural families of morphism $TTx \rightarrow x$ and $x \rightarrow Tx$ which give $T$ the monad structure. 

A warping  $(t, k)$ from $(X, B)$ to $T^\ast$ is an mw-monad on the category $B$ (see \cite{LS14}). The 2-cell $t$ corresponds to a family of functions 
\[T : B(x, Ty) \rightarrow B(Tx, Ty).\]
\noindent The 2-cell $k$ corresponds to specifying morphisms $K_x : x \rightarrow Tx$ of $B$, for each $x$. The warping axioms then translate to the equations
\[T(Tg\circ f) = Tg\circ Tf\]
\[TK_x = 1_{Tx}\]
\[f = Tf\circ K_x. \]

\noindent which are exactly the data and the axioms for an mw-monad. 

The construction of Theorem \ref{warptheo} is exactly the Kleisli construction for an mw-monad. The Theorem \ref{warptowreath} recaptures the one-to-one correspondence between monads and mw-monads.   
\end{example}

\section{Skew monads}\label{SkewMon}
Further we will work within the tricategory $\bc$. As usual, we will write as if $\bc$ were a gray category.

In \cite{LS12} the notion of a \textit{skew monoidale} within a monoidal bicategory was defined. A \textit{skew monad} in a tricategory $\bc$ is nothing else but a skew monoidale in the endo-hom monoidal bicategory $\bc(X, X)$, for an object $X$ in $\bc$. As the notion of skew monoidale is a laxification of the notion of monoid, so the notion of skew monad is a laxification of the notion of monad. More verbosely, a skew monad is similar to a monad except that the identities between 2-cells in the monad axioms are replaced by appropriately directed and suitably coherent structural non-invertible 3-cells. 

\begin{definition}
A (left) \textbf{skew monad} in $\bc$ is defined to be a skew monoidale in the endo-hom monoidal bicategory $\bc(X, X)$ of an object $X$ in $\bc$; besides the object $X$ it consists of an endomorphism $B : X \rightarrow X$, 2-cells $p : BB \Rightarrow B$ and $e : 1_X \Rightarrow B$, a 3-cell called an associator

\[\xymatrix{BBB \ar@{=>}[rr]^{Bp}\ar@{=>}[d]_{pB}&&BB\ar@{=>}[d]^{p}\\
BB\ar@{=>}[rr]_{p}="1"&&B \ar@3{>}^{\alpha}"1"+(0,11); "1"+(0,7) }
\]

\noindent and 3-cells called left and right unitors
\[
\xymatrix{&BB\ar@{=>}[rd]^{p}&\\
B\ar@{=>}[ur]^{Be}\ar@{=>}[rr]_{1_B}="1"&&B \ar@3{>}^{\lambda}"1"+(0,10); "1"+(0,6)
}
\qquad
\xymatrix{&BB\ar@{=>}[rd]^{p}&\\
B\ar@{=>}[ur]^{eB}\ar@{=>}[rr]_{1_B}="1"&&B \ar@3{>}^{\rho}"1"+(0,6); "1"+(0,10)
}
\]

\noindent satisfying the following five axioms:

\[\xymatrix{
p(pB)(BpB)\ar@3{>}[r]^{\alpha-}&p(Bp)(BpB)\ar@3{>}[r]^{-\alpha}& p(Bp)(BBp)\\
p(pB)(pBB) \ar@3{>}[u]^{-\alpha}\ar@3{>}[r]_{\alpha-}&p(Bp)(pBB)\ar@3{>}[r]_{\cong}&p(pB)(BBp) \ar@3{>}[u]_{\alpha-}\\
}\]

\[
\xymatrix{p(pB)(BeB)\ar@3{>}[r]^{\alpha-}&p(Bp)(BeB)\ar@3{>}[d]^{-\lambda}\\
p \ar@3{>}[u]^{-\rho}\ar@3{>}[r]_{1}&p}
\]

\[
\xymatrix{p(pB)(eBB)\ar@3{>}[r]^{-\lambda}\ar@3{>}[d]_{\alpha-}&p\\
p(Bp)(eBB)\ar@3{>}[r]_{\cong}&p(eB)p\ar@3{>}[u]_{\lambda-}}
\]

\[
\xymatrix{p\ar@3{>}[r]^<<<<<{-\rho}\ar@3{>}[d]_{\rho-}&p(Bp)(BBe)\\
p(Be)p\ar@3{>}[r]_{\cong}&p(pB)(BBe)\ar@3{>}[u]_{\alpha-}}
\]

\[
\xymatrix{p(Be)e\ar@3{>}[r]^{\cong}&\ar@3{>}[d]^{\lambda-}p(eB)e\\
e\ar@3{>}[u]^{\rho-}\ar@3{>}[r]_{1}&e}
\]

\end{definition}

\begin{example}\label{exmonoidale}
Every monoidal bicategory can be considered as a one-object tricategory, whereupon the latter is known as the suspension of the former. Suppose that $\bc$ is a suspension of a monoidal bicategory $(\cv, \otimes, I)$. Then, by definition, a skew monad in $\bc$ is the same as the skew monoidale in $\cv$.
\end{example}

\begin{example}\label{exskewmoncat}
A skew monoidale in the monoidal bicategory $(\Cat, \times, 1)$ is a skew monoidal category of \cite{Sz12}. It follows that a skew monad in the suspension of $\Cat$ is a skew monoidal category.
\end{example}

\begin{example}\label{exskewbicat}
Consider the tricategory of category-matrices $\Mat(\Cat)$ which is a higher dimensional analogue of the bicategory $\Span$ described in \cite{Ch15}. Its object are sets. The hom-bicategory $\Mat(\Cat)(X, Y)$ is defined to be the 2-category $[X\times Y, \Cat]$. Thus, a morphism $M : X \rightarrow Y$ is essentially a collection of small categories $M(x, y)$ indexed by the elements of the product $X \times Y$. The morphisms and 2-cells are given by indexed collections of functors and natural transformations. The horizontal composition is defined by the usual matrix multiplication formula
\[MN(x, y) = \coprod_zM(z, y)\times N(x, z).\]

\noindent The identity morphism $1_X$ is the collection of categories indexed by the elements of $X\times X$ which has the one-object category on the diagonal and the empty category everywhere else. Since its homs are strict 2-categories, $\Mat(\Cat)$ is in fact a $2\text{-}\Cat$-enriched bicategory.

A skew monad in $\Mat(\Cat)$ is a skew bicategory introduced in \cite{LS14}. Indeed, suppose that $(X, B, p, e, \alpha, \rho, \lambda)$ is a skew monad in $\Mat(\Cat)$. Then, $X$ is the set of objects of the skew bicategory.  For each pair of element $x$ and $y$ of $X$, $B(x, y)$ is the hom-category, i.e. its objects are morphisms $x \rightarrow y$ of the skew bicategory, and its morphisms are 2-cells between these in the skew bicategory. The 2-cell $p : AA \Rightarrow A$ amount to a family of functors
\[ A(z, y) \times A(x, z) \rightarrow A(x, y)\]

\noindent which sends a pair $(g, f)$ to their composite $g\circ f$. The 2-cell  $e : 1_X \Rightarrow A$ amounts to specifying an object $1_x$ in $A(x, x)$ for each $x$. The 3-cells $\alpha$, $\rho$ and $\lambda$ correspond to natural transformations whose components are
\[\xymatrix{(h\circ g)\circ f\ar[r]&h\circ (g\circ f)}\]
\[\xymatrix{1\circ f\ar[r]& f}\]
\[\xymatrix{f\ar[r]& f\circ 1}.\]

\noindent These are the non-invertible associator and unitor structural morphisms of the skew bicategory. The axioms of a skew monad translate to the axioms which these structural morphisms are required to satisfy.
\end{example}

\begin{example}
Let $\bc$ be the double suspension of a braided monoidal category $(\cb, \otimes, i)$. Then, a skew monad with $e = i$ is an \textit{augmented lax tricocycloid} (see \cite{LS12}). The 2-cell $p$ is the underlying object of the lax tricocycloid. The skew associator $\alpha$ determines the fusion operator $p\otimes p \rightarrow p\otimes p$. While, the skew unitors $\lambda$ and $\rho$ provide the \textit{augmentation} which consists of a counit $p \rightarrow i$ and a unit $i \rightarrow p$.
\end{example}

\section{Skew warpings}\label{SkewWarp}
Just as a skew monad is a laxification of a monad, so a skew warping is a laxification of a warping.

\begin{definition}
A \textbf{skew warping} from a monad $(X, B)$ to an endomorphism $A : X \rightarrow X$ consists of 2-cells $t : ABA \Rightarrow AB$ and $k : 1_X \Rightarrow AB$ and 3-cells

\[
\xymatrix@R=1em{
&ABBA\ar@{=>}[r]^{ApA}="1"&ABA\ar@{=>}[rd]^{t}& \\
ABABA \ar@{=>}[rd]_{ABt}\ar@{=>}[ru]^{tBA}&&&AB\\
&ABAB\ar@{=>}[r]_{tB}&ABB\ar@{=>}[ru]_{Ap}&\ar@3{>} "1"+(0,-10); "1"+(0,-15) ^\nu\\
}
\]

\[
\xymatrix{&ABA\ar@{=>}[rd]^{t}&\\
A\ar@{=>}[ur]^{kA}\ar@{=>}[rr]_{Ae}="1"&&AB \ar@3{>}^{\nu_0}"1"+(0,10); "1"+(0,6)
}
\qquad\qquad
\xymatrix{ABAB \ar@{=>}[r]^{tB}&ABB\ar@{=>}[d]^{Ap}\\
AB\ar@{=>}[u]^{ABk}\ar@{=>}[r]_{1_{AB}}="1"&AB \ar@3{>}_{\kappa}"1"+(0,8); "1"+(0,12) }
\]

\noindent satisfying the following axioms

\[
\xymatrix{t(ApA)(tBA)(ApABA)(tBABA)\ar@3{>}[d]_{\nu-}\ar@3{>}[r]^{-\nu}&t(ApA)(ApBA)(tBBA)(ABtBA)\ar@3{>}[d]^{-\alpha-}\\
(Ap)(tB)(ABt)(ApABA)(tBABA)\ar@3{>}[d]_{\cong}&t(ApA)(ABpA)(tBBA)(ABtBA)\ar@3{>}[d]^{\cong}\\
(Ap)(tB)(ApAB)(tBAB)(BABt)\ar@3{>}[d]_{-\nu-}&t(ApA)(tBA)(ABApA)(ABtBA)\ar@3{>}[d]^{\nu-}\\
(Ap)(ApB)(tBB)(ABtB)(ABABt)\ar@3{>}[d]_{\alpha-}&(Ap)(tB)(ABt)(ABApA)(ABtBA)\ar@3{>}[d]^{-\nu}\\
(Ap)(ABp)(tBB)(ABtB)(ABABt)&(Ap)(tB)(ABAp)(ABtB)(ABABt)\ar@3{>}[l]^{\cong}\\
}
\]

\[
\xymatrix{(Ap)(tB)(ABt)(ABkA)\ar@3{>}[r]^{-\nu_0}&(Ap)(tB)(ABAe)\ar@3{>}[dd]^{\cong}\\
t(ApA)(tBA)(ABkA)\ar@3{>}[u]^{\nu-}&\\
t\ar@3{>}[u]^{\kappa-}\ar@3{>}[r]_{-\nu}&(Ap)(ABe)t\\}
\]

\[
\xymatrix{t(ApA)(tBA)(kABA)\ar@3{>}[d]_{-\nu_0}\ar@3{>}[r]^{\nu-}&(Ap)(tB)(ABt)(kABA)\ar@3{>}[d]^{\cong}\\
t(ApA)(AeBA)\ar@3{>}[d]_{-\lambda}&(Ap)(tB)(kAB)t\ar@3{>}[d]^{-\nu_0-}\\
t&(Ap)(AeB)t\ar@3{>}[l]^{\lambda-}\\
}
\]

\[
\xymatrix{(Ap)(tB)(ApAB)(tBAB)(ABABk)\ar@3{>}[r]^{-\nu-}&(Ap)(ApA)(tBB)(ABtA)(ABABk)\ar@3{>}[d]^{\alpha-}\\
(Ap)(tB)(ABk)(Ap)(tB)\ar@3{>}[u]^{\cong}&(Ap)(ABp)(tBB)(ABtA)(ABABk)\ar@3{>}[d]^{\cong}\\
(Ap)(tB)\ar@3{>}[u]^{\kappa-}\ar@3{>}[r]_{\kappa-}&(Ap)(tB)(ABp)(ABtA)(ABABk)\\}
\]

\[
\xymatrix{(Ap)(tB)(ABk)k\ar@3{>}[r]^{\nu_0-}&(Ap)(AeB)k\ar@3{>}[dd]^{\lambda-}\\
(Ap)(tB)(kAB)k\ar@3{>}[u]^{\cong}&\\
k\ar@3{>}[u]^{\kappa-}\ar@3{>}[r]_{}&k\\}
\]

\end{definition}

\begin{theorem}\label{skewwarptheo}
Given a skew warping $(t, k, \nu, \nu_0, \kappa)$ from a monad $(X, B)$ to an endomorphism $A$, there is a skew monad structure on the composite $AB$ with the skew unit $k : 1_X \Rightarrow AB$, the skew multiplication $(Ap)(tB) : ABAB \Rightarrow AB$, and the skew associator and the skew unitors defined from $\nu, \nu_0$ and $\kappa$ in the straightforward way. 
\end{theorem}

\begin{example}
Let $\bc$ be the suspension of a monoidal bicategory $(\cv, \otimes, I)$. By Example \ref{exmonoidale} a skew monad $(X, B)$ in $\bc$ is a skew monoidale in $\cv$. A warping in our sense from $(X, B)$ to $I$ is essentially the same as the skew warping  in the sense of \cite{LS12} on the corresponding monoidale in $\cv$. 
\end{example}

\begin{example}
Specializing the previous example by taking $\cv$ to be the monoidal 2-category $\Cat$, we recapture the notion of warping on a skew monoidal category \cite{LS12}.
\end{example}

\begin{example}
Take $\bc$ to be the bicategory $\Mat(\Cat)$. We have a pseudofunctor $\Set \rightarrow \Mat(\Cat)$ which sends a map $F : X \rightarrow Y$ to a category matrix $F^\ast$ for which $F^\ast(x, y)$ is the one-object category whenever $Fx = y$ and the empty category otherwise. $F^\ast$ has a dually defined left adjoint $F_\ast$. We have 
\[F^\ast M(z, x) = M(z, Fx)\]
\[NF^\ast(x, z) = N(Fx, z),\]

\noindent from which it follows that, a 2-cell  $MF^\ast \rightarrow F^\ast M$ amounts to a collection of functors $M(x, y) \rightarrow M(Fx, Fy)$.
 
By Example \ref{exskewbicat}, a skew monad $(X, B)$ in $\Mat(\Cat)$ is a skew bicategory, for which we also write $B$. 

A skew warping $(t, k, \nu, \nu_0, \kappa)$ from $(X, B)$ to $T^\ast : X \rightarrow X$ is the same as a skew warping on the skew bicategory $B$ in the sense of \cite{LS14}. $T$ is a map of sets $X \rightarrow X$. The 2-cell $t$ amounts to a family of functors
\[T : B(x, Ty) \rightarrow B(Tx, Ty),\]

\noindent The 2-cell $k$ amounts to specifying a family of objects $K_x$ in $B(x, Tx)$. The 3-cells $\nu$, $\nu_0$ and $\kappa$ translate to natural transformations  with components
\[\xymatrix{T(Tg\circ f)\ar[r]&Tg\circ Tf}\]
\[\xymatrix{TK_x\ar[r]&1_{Tx}}\]
\[\xymatrix{f\ar[r]&Tf\circ K_x},\]
 
\noindent which satisfying the axioms of a skew warping on a skew bicategory.

The construction of Theorem \ref{skewwarptheo} recaptures the Kleisli construction for a skew warping on a skew bicategory described in Section 5 in \cite{LS14}.  
\end{example}

\begin{example}
Let $\bc$ be a double suspension of a braided monoidal category $(\cb, \otimes, i)$. Consider a trivial skew monad, i.e. the one all of whose data consists of identity cells. A warping from the trivial monad to the only morphism of $\bc$ with $k = i$ and $j = i$ is essentially the same as an augmented lax tricocycloid with the underlying object $q$.
\end{example}

\section{Algebras}\label{algebras}
First we consider the bicategorical case. There exists a natural notion of an actee for a warping. 

\begin{definition}\label{algebra}
An algebra of a warping $(t, k)$ from a monad $(X, B)$ to an endofunctor $A$ in a bicategory $\bb$, consists of an object $Y$ of $B$, a morphism $M : X \rightarrow Y$ and a 2-cell $m : MBA \Rightarrow MB$ satisfying the following equalities

\[
\xymatrix@R=1em{
&MBBA\ar@{=>}[r]^{MpA}&MBA\ar@{=>}[dr]^{m}&\\
MBABA \ar@{=>}[ru]^{MBt}\ar@{=>}[dr]_{mBA}&&&MB\\
&MBAB\ar@{=>}[r]_{mB}&MBB\ar@{=>}[ur]_{Mp}&\\
}
\]

\[\xymatrix{MBAB \ar@{=>}[r]^{mB}&MBB\ar@{=>}[d]^{Mp}\\
MB\ar@{=>}[u]^{MBk}\ar@{=>}[r]_{1_{AB}}&MB}
\]
\end{definition}

\begin{example}
Suppose that $(t, k)$ is a warping from a monad $(X, B)$ to an endomorphism $T^\ast$ within a bicategory $\Span$. In other words, by Example \ref{exwarp}, it is an mw-monad on the category $B$. Suppose that $Y$ is the one-element set, and $a$ is the map $Y \rightarrow X$ corresponding to choosing an object $a$ of $B$. A 2-cell $m : a^\ast BT^\ast \rightarrow a^\ast B$ amounts to a family of maps $E : B(z, Ta) \rightarrow B(z, a)$. The algebra axioms become the identities
\[E(Eg\circ f) = Tg\circ Ef\]
\[g = Eg\circ K_x\]

\noindent where $g: y \rightarrow a$ and $f : x \rightarrow Ty$. This means that an algebra for the warping is the same as an algebra for the mw-monad $T$. The latter itself is the same as an algebra for the corresponding usual monad. The usual monad algebra structure is given by the morphism $E1_{Tx} : Tx \rightarrow x$. 
\end{example}

Now we switch to the tricategorical setting. An algebra of a skew warping is a laxification of an algebra of a warping. 

\begin{definition}
An algebra of a skew warping $(t, k, \nu, \nu_0, \kappa)$ from a skew monad $(X, B)$ to an endofunctor $A$ in a tricategory $\bc$, consists of an object $Y$ of $B$, a morphism $M : X \rightarrow Y$, a 2-cell $m : MBA \Rightarrow MB$ and 3-cells

\[
\xymatrix@R=1em{
&MBBA\ar@{=>}[r]^{MpA}="1"&MBA\ar@{=>}[dr]^{m}&\\
MBABA \ar@{=>}[rd]_{MBt}\ar@{=>}[ur]^{mBA}&&&MB\\
&MBAB\ar@{=>}[r]_{mB}&MBB\ar@{=>}[ur]_{Mp}& \ar@3{>} "1"+(0,-10); "1"+(0,-15) ^\nu\\
}
\]

\[\xymatrix{MBAB \ar@{=>}[r]^{mB}&MBB\ar@{=>}[d]^{Mp}\\
MB\ar@{=>}[u]^{MBk}\ar@{=>}[r]_{1_{MB}}="1"&MB \ar@3{>}_{\kappa}"1"+(0,8); "1"+(0,12)}
\]

\noindent coherent in the suitable way.
\end{definition}

\begin{example}
Suppose that $(X, B)$ is a skew monad in the tricategory $\Mat(\Cat)$, and $(t, k)$ is a skew warping from it to an endofunctor $T^\ast$. In other words it is a skew warping on the skew bicategory $B$. Suppose that $Y$ is the one-element set, and $a$ is the map $Y \rightarrow X$ corresponding to choosing an object $a$ of $B$. A 2-cell $m : a^\ast BT^\ast \rightarrow a^\ast B$ amounts to a family of functors $E : B(z, Ta) \rightarrow B(Tz, Ta)$. The structural 3-cells of the algebra amount to morphisms
\[\xymatrix{E(Eg\circ f) \ar[r]& Tg\circ Ef}\]
\[\xymatrix{g \ar[r]& Eg\circ K_x}\]

\noindent where $g: y \rightarrow a$ and $f : x \rightarrow Ty$, satisfying exactly the axioms of an algebra for a skew warping on a skew bicategory in the sense of \cite{LS14}.
\end{example}

%\bibliographystyle{alpha}
%\bibliography{bib}
\end{document}